\documentclass[12pt]{article}

\usepackage{amssymb,latexsym,amsfonts,amsmath,amsthm,verbatim,eucal,hyperref}

\theoremstyle{plain}
\newtheorem*{thm*}{Theorem}
\newtheorem{clm}{Claim}
\theoremstyle{remark}
\newtheorem*{rmk*}{Remark}
\def\P{\mathbb{P}}
\def\E{\mathbb{E}}
\def\R{\mathbb{R}}
\def\S{S^{n-1}}
\def\eps{\varepsilon}

\begin{document}

\title{An extension of a Bourgain--Lindenstrauss--Milman inequality}
\author{Omer Friedland, Sasha Sodin}
\maketitle

\begin{abstract}
\noindent Let $\| \cdot \|$ be a norm on $\R^n$. Averaging $\| (\eps_1 x_1, \cdots, \eps_n x_n) \|$
over all the $2^n$ choices of $\overrightarrow{\eps} = (\eps_1, \cdots, \eps_n) \in \{ -1, +1 \}^n$,
we obtain an expression $\|| x \||$ which is an unconditional norm on $\R^n$.

\vspace{2mm}
\noindent Bourgain, Lindenstrauss and Milman \cite{blm} showed that, for a certain (large) constant
$\eta > 1$, one may average over $\eta n$ (random) choices of $\overrightarrow{\eps}$ and obtain
a norm that is isomorphic to $\|| \cdot \||$. We show that this is the case for any $\eta > 1$.
\end{abstract}

\section{Introduction}
\footnotetext[1]{[omerfrie; sodinale]@post.tau.ac.il; address: 
School of Mathematical Sciences, Tel Aviv University, Ramat Aviv, 
Tel Aviv 69978, Israel}

Let $(E, \| \cdot \|)$ be a normed space, and let $v_1, \cdots, v_n \in E \setminus \{ 0 \}$.
Define a norm $\|| \cdot \||$ on $\R^n$:
\begin{equation}\label{1}
\|| x \|| = \E \, \| \sum \eps_i x_i v_i \|~,
\end{equation}
where the expectation is over the choice of $n$ independent random
signs $\eps_1, \cdots, \eps_n$.  This is an {\em unconditional} norm; that is,
\[ \|| (x_1, x_2, \cdots, x_n) \|| = \|| (|x_1|, |x_2|, \cdots, |x_n|) \||~.\]

The following theorem states that it is sufficient to average $O(n)$, rather than $2^n$,
terms in (\ref{1}), in order to obtain a norm that is isomorphic to $\|| \cdot \||$ (and in
particular approximately unconditional).

\begin{thm*}
Let $N = (1+\xi)n$, $\xi > 0$, and let
\[ \{ \eps_{ij} \, \big| \, 1 \leq i \leq n, 1 \leq j \leq N \} \]
be a collection of independent random signs. Then
\begin{equation*}
\P \left\{ \forall x \in \R^n \,\,
    c(\xi) \, |||x|||
        \le \frac{1}{N} \sum_{j=1}^N\|\sum_{i=1}^n \eps_{ij} x_i v_i\|
        \le C(\xi) |||x||| \right\}
    \ge 1-e^{-c' \xi n}~,
\end{equation*}
where
\[ c(\xi) =
    \begin{cases}
        c \xi^2,        & 0 < \xi < 1 \\
        c,              &1 \leq \xi < C'' \\
        1 - C'/\xi^2,   &C'' \leq \xi
    \end{cases}, \quad
C(\xi) =
    \begin{cases}
        C,              &0 < \xi < C'' \\
        1 + C''/\xi^2,  &C'' \leq \xi
    \end{cases}~,\]
and $c,c',C,C',C''>0$ are universal constants (such that $1 - C'/C''^2 \geq c$,
$1 + C'/{C''}^2 \leq C$).
\end{thm*}

This extends a result due to Bourgain, Lindenstrauss and Milman \cite{blm},
who considered the case of large $\xi$ ($\xi \geq C''$); their proof
makes use of the Kahane--Khinchin inequality. Their argument yields the upper
bound for the full range of $\xi$, so the innovation is in the lower bound
for small $\xi$.

With the stated dependence on $\xi$, the corresponding result for the scalar case
$\dim E = 1$ was proved by Rudelson \cite{r}, improving previous bounds
on $c(\xi)$ in \cite{lprt,afm,afms}; see below. This is one of the two main ingredients
of our proof, the second one being Talagrand's concentration inequality \cite{t}
(which, as shown by Talagrand, also implies the Kahane--Khinchin inequality).

\vspace{2mm}
{\em Acknowledgement}: We thank our supervisor Vitali Milman for his support
and useful discussions.

\section{Proof of Theorem}

Let us focus on the case $\xi < 1$; the same method works (in fact, in
a simpler way) for $\xi \geq 1$.

Denote $\|| x \||_N = \frac{1}{N} \sum_{j=1}^N\|\sum_{i=1}^n \eps_{ij} x_i v_i\|$;
this is a random norm depending on the choice of $\eps_{ij}$. Let
$\S_{|||\cdot|||}=\{x\in\R^n:|||x|||=1\}$ be the unit sphere of $(\R^n, \|| \cdot \||)$;
we estimate
\begin{equation}\label{unb}\begin{split}
&\P\left\{ \forall x\in\S_{|||\cdot|||}, \, \,c \xi^2 \le|||x|||_N \le C \right\} \\
&\quad\geq 1 - \P \left\{ \exists x\in\S_{||| \cdot |||}, \, |||x|||_N>C \right\} \\
&\qquad- \P \left\{ \left(\forall y\in\S_{||| \cdot |||}, \, |||y|||_N\le C \right)
        \wedge \left(\exists x\in\S_{||| \cdot |||} , |||x|||_N<c\xi^2 \right) \right\}~.
\end{split}\end{equation}

\vskip 10pt \noindent \textbf{Upper bound:} Let us estimate the first term
\[ \P\left\{ ~\exists x\in\S_{||| \cdot |||} ,|||x|||_N>C \right\}~. \]

\begin{rmk*}
As we mentioned, the needed estimate follows from the argument in \cite{blm};
for completeness, we reproduce a proof in the similar spirit.
\end{rmk*}

\begin{thm*}[Talagrand \cite{t}]\label{talagrand}
Let $w_1, \cdots, w_n \in E$ be vectors in a normed space $(E, \| \cdot \|)$,
and let $\eps_1,\cdots,\eps_n$ be independent random signs. Then for any $t > 0$
\begin{equation}\label{teq}
\P \left\{ \left| \|\sum_{i=1}^n \eps_i w_i\|
    - \E \, \|\sum_{i=1}^n \eps_i w_i\| \right| \ge t \right\}
    \le C_1 e^{-c_1 t^2/\sigma^2},
\end{equation}
where $c_1,C_1>0$ are universal constants, and
\[ \sigma^2
    = \sigma^2(w_1, \cdots, w_n)
    =\sup\left\{\sum_{i=1}^n\varphi(w_i)^2 \, \big{|} \,
    \varphi \in E^*, \, \|\varphi\|^*\le1\right\}~.\]
\end{thm*}

\begin{rmk*}
Talagrand has proved (\ref{teq}) with the median
$\operatorname{Med} \|\sum_{i=1}^n \eps_i w_i\| $ rather than the
expectation; one can however replace the median by the expectation
according to the proposition in Milman and Schechtman
\cite[Appendix V]{ms}.
\end{rmk*}

For $x = (x_1, \cdots, x_n) \in \R^n$, denote
\[ \sigma^2(x) = \sigma^2(x_1 v_1, \cdots, x_n v_n)~.\]

\begin{clm} $\sigma$ is a norm on $\R^n$ and
$\sigma(x) \leq C_2 \||x\||$ for any $x \in \R^n$.
\end{clm}

\begin{proof} The first statement is trivial. For the second one,
note that
\[ \||x\||
    = \E \, \| \sum \eps_i x_i v_i \|
    \geq \E \left|\varphi(\sum \eps_i x_i v_i)\right|
    = \E \left|\sum \eps_i \varphi(x_i v_i)\right|,
    \quad \| \varphi \|^\ast \leq 1~. \]
Now, by the classical Khinchin inequality,
\begin{equation}\label{khin}
\sqrt{\sum y_i^2}
    \geq \E \left|\sum \eps_i y_i \right|
    \geq C_2^{-1} \sqrt{\sum y_i^2} \end{equation}
(see Szarek \cite{s} for the optimal constant $C_2 = \sqrt{2}$).
Therefore
\[ \|| x \||
    \geq C_2^{-1} \sup_{\|\varphi\|^\ast \leq 1} \sqrt{\sum \varphi(x_i v_i)^2}
    = C_2^{-1} \sigma(x)~.\]
\end{proof}

By the claim and Talagrand's inequality, for every (fixed)
$x \in S^{n-1}_{\|| \cdot \||}$
\[ \P \left\{  \|\sum_{i=1}^n \eps_i x_i\| \geq t \right\} \leq C_1 \exp(-c_2 t^2)~.\]
Together with a standard argument (based on the exponential Chebyshev
inequality), this implies (for $t$ large enough):
\[ \P \left\{  \frac{1}{N} \sum_{j=1}^N \|\sum_{i=1}^n \eps_{ij} x_i\| \geq t \right\}
    \leq \exp(-c_3 t^2 N)~.\]
In particular, for $t = C_3 \geq \sqrt{4/c_3}$ the left-hand side is smaller than
$12^{-N} < 6^{-n} 2^{-N}$.

The following fact is well-known, and follows for example from
volume estimates (cf. \cite{ms}).
\begin{clm} For any $\theta > 0$, there exists a $\theta$-net
$\mathcal{N}_\theta$ with respect to $\||\cdot\||$ on $S^{n-1}_{\||
\cdot \||}$ of cardinality $\# \mathcal{N}_\theta \leq
(3/\theta)^n$.
\end{clm}

For now we only use this for $\theta=1/2$. By the above,
with probability greater than $1 - 2^{-N}$, we have: $\|| x \||_N \leq C_3$
simultaneously for all $x \in \mathcal{N}_{1/2}$.

Representing an arbitrary unit vector $x \in S^{n-1}_{\|| \cdot
\||}$ as
\[ x = \sum_{k=1}^\infty a_k x^{(k)},
    \quad |a_k| \leq 1/2^{k-1}, \, x^{(k)} \in \mathcal{N}_{1/2}~, \]
we deduce: $\|| x \||_N \leq 2 C_3$, and hence finally:
\begin{equation}\label{u}
\P \left\{ \exists x\in\S_{||| \cdot |||}, \, |||x|||_N>C
\right\}
    \leq 2^{-N}
\end{equation}
(for $C = 2C_3$).

\vskip 10pt \noindent \textbf{Lower bound:} Now we turn to the second term
\begin{equation*}\label{lower}
\P \left\{ \left( \forall y\in\S_{||| \cdot |||} ,|||y|||_N \le C
\right)
    \wedge \left(\exists x\in\S_{||| \cdot |||} ,|||x|||_N<c\xi^2\right) \right\}~.
\end{equation*}

For $\sigma_0$ (that we choose later), let us decompose
$ S^{n-1}_{\|| \cdot \||} = U \uplus V$, where
\[ U = \left\{ x \in S^{n-1}_{\|| \cdot \||} \, \Big| \,
        \sigma(x) \geq \sigma_0 \right\},
    \quad V = \left\{ x \in S^{n-1}_{\|| \cdot \||} \, \Big| \,
        \sigma(x) < \sigma_0 \right\}~. \]

\vskip 10pt \noindent  Recall the following result (mentioned in the introduction);
we use the lower bound that is due to Rudelson \cite{r}.
\begin{thm*}[\cite{lprt,afm,afms,r}]\label{rudelson}
Let $N = (1+\xi)n$, $0 < \xi < 1$, and let
\[ \{ \eps_{ij} \, \big| \, 1 \leq i \leq n, 1 \leq j \leq N \} \]
be a collection of independent random signs. Then
\begin{equation*}
\P \left\{ \forall y \in \R^n \,\,
    c_4 \xi^2 \, |y|
        \le \frac{1}{N} \sum_{j=1}^N |\sum_{i=1}^n \eps_{ij} y_i|
        \le C_4 |y| \right\}
    \ge 1-e^{-c_4' \xi n}~,
\end{equation*}
where $c_4,c_4',C_4>0$ are universal constants, and $|\cdot|$ is the standard Euclidean norm.
\end{thm*}

\begin{rmk*} By the Khinchin inequality (\ref{khin}), this is indeed the scalar case of Theorem~1
for $0 < \xi < 1$.
\end{rmk*}

Thence with probability $\geq 1 - e^{-c_4' \xi n}$ the following
inequality holds for all $x\in U$ (simultaneously):
\begin{equation}\label{l-U}\begin{split}
|||x|||_N
    &\ge \frac{1}{N}\sum_{j=1}^N\left|\varphi(\sum_{i=1}^n\eps_{ij} x_i v_i)\right|
    = \frac{1}{N}\sum_{j=1}^N \left|\sum_{i=1}^n \eps_{ij} x_i\varphi(v_i)\right| \\
    &\ge c_4\xi^2\sigma(x) \geq c_4 \xi^2 \sigma_0~.
\end{split}\end{equation}

\vskip 10pt \noindent Now let us deal with vectors $x\in V$. Let
$\mathcal{N}_\theta$ be a $\theta$-net on $S^{n-1}_{\|| \cdot \||}$
(where $\theta$ will be also chosen later). For $x' \in
\mathcal{N}_\theta$ such that $\|| x - x' \||\leq \theta$,
$\sigma(x') \leq \sigma_0 + C_2\theta$ by Claim 1. Therefore by
Talagrand's inequality (\ref{teq}),
\[ \P \left\{ \|\sum_{i=1}^n \eps_i x'_i v_i\| < 1/2 \right\}
    \leq C_1 \exp(-c_1/(4 (\sigma_0 + C_2\theta)^2))~,\]
and hence definitively
\begin{multline*}
\P \left\{ \frac{1}{N} \sum_{j=1}^N \|\sum_{i=1}^n \eps_{ij} x'_i v_i\| < 1/4 \right\}
    \leq 2^N \left\{ C_1 \exp\left(-\frac{c_1}{4 (\sigma_0+C_2\theta)^2}\right) \right\}^{N/2} \\
    = \exp \left\{ - \left(\frac{c_1}{8(\sigma_0+C_2\theta)^2} - \log (2\sqrt{C_1}) \right) N \right\}~.
\end{multline*}
Let $\sigma_0 = C_2 \theta$, and choose $0 < \theta < 1/(8C)$ so that
\[ \frac{c_1}{32 C_2^2 \theta^2} - \log (2 \sqrt{C_1}) > \log 2 + \log (3/\theta)~.\]
Then the probability above is not greater than $2^{-N} (\theta/3)^N
< 2^{-N} / \# \mathcal{N}_\theta$ (by Claim~2). Therefore with
probability $\geq 1 - 2^{-N}$ we have:
\[ \|| x' \||_N \geq 1/4 \quad \text{for $x' \in \mathcal{N}_\theta$ such that $\||x - x'\|| < \theta$
for some $x \in V$}~.\]
Using the upper bound (\ref{u}), we infer:
\begin{equation}\label{l-V}\begin{split}
\||x\||_N
    &\geq \||x'\||_N - \||x' - x\||_N \\
    &\geq 1/4 - C/8C = 1/4 - 1/8 = 1/8~, \quad x \in V~.
\end{split}\end{equation}

The juxtaposition of (\ref{unb}), (\ref{u}), (\ref{l-U}), and
(\ref{l-V}) concludes the proof. $\hfill\square$

\end{document}